%% file: main.tex
\documentclass[conference,a4paper]{APSIPA2021}
\usepackage{amsmath,graphicx,amssymb}
\usepackage{calc,bm,url,color,theorem,graphicx,cite,shortcuts_FW}
\usepackage{psfrag,subfigure,float,hyperref}
\usepackage{paralist}
\usepackage{mathtools}

\usepackage[flushleft]{threeparttable}
\usepackage{microtype} 
\usepackage{nicefrac} 
\usepackage{diagbox}
\usepackage{wrapfig}

\usepackage{multirow, booktabs}
\usepackage{algorithm}
\usepackage{algorithmic}

\usepackage{etoolbox}
\makeatletter
\patchcmd{\@makecaption}
  {\scshape}
  {}
  {}
  {}
\makeatletter
\patchcmd{\@makecaption}
  {\\}
  {.\ }
  {}
  {}
\makeatother
\def\tablename{Table}

\newtheorem{Theorem}{Theorem}
\newtheorem{Assumption}{Assumption}

\usepackage[utf8]{inputenc}
\usepackage{pgfplots}
\DeclareUnicodeCharacter{2212}{−}
\usepgfplotslibrary{groupplots,dateplot}
\usetikzlibrary{patterns,shapes.arrows}
\pgfplotsset{compat=newest}

\title{An Empirical Study on Compressed Decentralized Stochastic Gradient Algorithms with Overparameterized Models}

\author{%
\authorblockN{%
Arjun Ashok Rao\authorrefmark{1} and
Hoi-To Wai\authorrefmark{1}
}
\authorblockA{%
\authorrefmark{1}
Dept.~of SEEM; The Chinese University of Hong Kong; Shatin, Hong Kong\\
E-mail: \texttt{arjunrao@link.cuhk.edu.hk, htwai@cuhk.edu.hk}
}
}

\thanks{SSS}
\begin{document}
%
\maketitle
\begin{abstract}
 This paper considers decentralized optimization with application to machine learning on graphs. The growing size of neural network (NN) models has motivated prior works on decentralized stochastic gradient algorithms to incorporate communication compression. On the other hand, recent works have demonstrated the favorable convergence and generalization properties of overparameterized NNs. In this work, we present an empirical analysis on the performance of compressed decentralized stochastic gradient (DSG) algorithms with overparameterized NNs. Through simulations on an MPI network environment, we observe that the convergence rates of popular compressed DSG algorithms are robust to the size of NNs. Our findings suggest a gap between theories and practice of the compressed DSG algorithms in the existing literature.
\end{abstract}
\begin{keywords}
decentralized optimization, communication efficiency, overparameterized models
\end{keywords}
\section{Introduction}
\label{sec:intro}
A popular approach for enabling scalable machine learning is applying decentralized algorithms to tackle the large-scale optimization problem through collaboration between a group of workers connected on a network/graph. 
Let $N \in \NN$ be the number of workers and $d \in \NN$ be the problem's dimension, we consider the following unconstrained optimization problem:
\beq \label{eq:opt_int}
\min_{ \prm \in \RR^d }~J(\prm) \quad \text{where} \quad J(\prm) := \frac{1}{N}\sum_{i=1}^N J_i( \prm ),
\eeq
where $J_i : \RR^d \rightarrow \RR$ is a continuously differentiable function representing the private data held by the $i$th worker. Furthermore, the $N$ workers are connected on an undirected graph denoted by $G = (V,E)$, where $V = [N] = \{1,...,N\}$ is the set of workers and $E \subseteq V \times V$ is the set of edges of $G$ with self loops such that $(i,i) \in E$ for all $i \in V$. Let $G$ be a connected graph, we note that \eqref{eq:opt_int} is equivalent to the consensus optimization problem:
\beq
\label{eq:opt_con}
\min_{ \prm_i \in \RR^d, i \in V} ~\sum_{i=1}^N J_i( \prm_i ) ~~\text{s.t.}~~\prm_i = \prm_j,~\forall~(i,j) \in E,
\eeq
where $\prm_i \in \RR^d$ is a private/local variable held by the $i$th worker. The applications of \eqref{eq:opt_int} include decentralized regression, sensor fusion for wireless sensor network \cite{sayed2014adaptation}, etc.

In this paper, we are concerned with the application of \eqref{eq:opt_int} to machine learning (ML) tasks via training a neural network (NN) model. Following the design of \eqref{eq:opt_int}, our goal is to train a common model $\prm$ at all workers. For example, if we consider a supervised learning problem for classification, the $i$th private function takes the form of the empirical risk: 
\beq \label{eq:ji}
J_i( \prm ) = \frac{1}{ |{\cal D}_i| } \sum_{j=1}^{|{\cal D}_i|} {\sf loss}( f( {\bm x}_j; \prm) ; y_j ) ,
\eeq
where ${\bm x}_j \in \RR^f$ and $y_j \in \RR$ are the $j$th feature and label known by worker $i$, respectively, and $|{\cal D}_i|$ is the number of samples held by worker $i$. The loss function ${\sf loss}(\cdot)$ can be taken as the cross-entropy, or the quadratic loss. The nonlinear function $f({\bm x} ; \prm)$ is the output of a neural network, e.g., a two-layer neural network with ReLU activation is given by 
\beq \label{eq:fi}
f({\bm x} ;\prm) = \frac{ 1 }{ \sqrt{m} } \sum_{j=1}^m b_j \max\{ 0,  \pscal{ {\bm x} }{ \prm^{(j)} } \},
\eeq
where $b_j$ is the $j$th output weight and we have defined the parameters as $\prm = ( \prm^{(1)}, ..., \prm^{(m)} ) \in \RR^{mf}$ such that $d= mf$. 
Notice that despite its simplicity, the NN architecture \eqref{eq:fi} exhibits good representation power provided that $m \rightarrow \infty$ \cite{bach2017breaking}. 

To tackle \eqref{eq:opt_int} when only local communications are allowed, decentralized first-order optimization algorithms have been developed for over a decade \cite{nedic2009distributed,sayed2014adaptation,shi2015extra}. The main idea behind these algorithms is a simultaneous ``consensus + optimize'' strategy where workers communicate with each other to reach a common model $\prm$ (i.e., achieving consensus) while optimizing their local models $\prm_i$ via gradient steps on the private functions. For a non-convex optimization model, the convergence of decentralized optimization algorithms has been analyzed in \cite{bianchi2012convergence,di2016next,hong2017prox} under the deterministic setting. Recent works have also analyzed the extension to stochastic gradient (SG) based methods \cite{lian2017can,zeng2018nonconvex,tang2018d,sun2020improving,xin2021improved}; also see \cite{chang2020distributed} for a recent overview. 

A major obstacle in applying the above algorithms to ML tasks \eqref{eq:ji} lies in the \emph{communication overhead}. Taking the decentralized gradient (DGD) method from \cite{nedic2009distributed} as an example, each training iteration requires workers to transmit the entire $d$-dimensional local model. Notice that the state-of-the-art NN models are typically large with $d \gg 1$, e.g., the VGG16 NN comprises of $d \approx 1.38 \times 10^8$ parameters \cite{simonyan2014very}. When applying the DGD method to train such a model, workers would be required to send 526 MB of data over the network \emph{per-iteration}; This may cause a significant slowdown to the practical convergence of the algorithms.

In light of the above, previous works have been motivated to develop \emph{communication efficient} variants of decentralized SG algorithms. To list a few examples, \cite{koloskova2019decentralized,Koloskova*2020Decentralized} proposed the CHOCO-SGD algorithm which compresses incremental vectors at the workers prior to transmission; also see its generalization in \cite{kovalev2021linearly}; \cite{tang2018communication} proposed the ECD-PSGD algorithm based on an extrapolation technique. It is worth noting that the above algorithms share similarities with the compression techniques proposed for distributed SG such as \cite{alistarh2017qsgd}.

On the other hand, recent works have studied the \emph{global convergence} of centralized SG algorithms for training \emph{overparameterized} NN models \cite{jacot2018neural} via the non-convex optimization \eqref{eq:opt_int}, \eqref{eq:ji}. Particularly, these works consider simple architectures such as \eqref{eq:fi} and show that the optimal weights $\prm^\star$ to \eqref{eq:opt_int} can be found using SG algorithms when $m \rightarrow \infty$, e.g., \cite{arora2019fine, bakshi2019learning}. Interestingly, it is demonstrated that the convergence rate for the $L^2$ function distance will be independent of $m$, implying that there is no loss in iteration complexity (for centralized learning) when deploying \emph{overparameterized} NNs. 

The above observation illustrates a dilemma when deciding the NN architecture to be used for decentralized training: on one hand, it is desired to adopt a compact NN model to reduce communication cost, on the other hand, using an \emph{overparameterized} NN model provides enticing theoretical convergence properties and higher representation power. 
The aim of this paper is to empirically study the effects of overparameterization on the performance of (communication efficient) decentralized SG algorithms. Particularly, we verify the claim that \emph{decentralized SG algorithms are practical (in terms of overall communication efficiency), even for training overparameterized NNs}. Our contributions are as follows:
\begin{itemize}
\item In Section~\ref{sec:alg}, we highlight the pitfalls of existing theories for training NNs with communication-efficient decentralized SG algorithms by showing that the curse of dimensionality persists. We argue that the latter can be an artifact of the analysis as it does not exploit the structure of the NN training problems. 
\item To support our claim, in Section~\ref{sec:num}, we concentrate on the CHOCO-SGD algorithm \cite{Koloskova*2020Decentralized} and experiment with this algorithm in various setting, providing the \emph{first set of experiments} that highlights on the effects of an increasing NN width $m$ in the model \eqref{eq:fi}. We show that, when the amount of information transmitted per iteration is fixed, the algorithm achieves a lower training loss / testing error at the same speed (or even faster) as the NN's width is increased.  
\end{itemize}
In addition, we discuss the intuitive reason behind the current gap in the theories and suggest potential solutions to fix it. 

\section{Decentralized Learning: Algorithms and Communication Efficiency} \label{sec:alg}
In this section, we discuss the algorithms for decentralized learning/training and their communication efficient variants. Furthermore, we will review the convergence guarantees of these algorithms and discuss how they may fail to yield meaningful insights for training overparameterized NNs. 

To fix ideas, we define the doubly stochastic mixing matrix ${\bm W} \in \RR_+^{N \times N}$ satisfying the row/column sum condition ${\bm W} {\bf 1} = {\bm W}^\top {\bf 1} = {\bf 1}$; it respects the graph topology such that $W_{ij} = W_{ji} = 0$ whenever $(i,j) \notin E$; moreover, it satisfies the fast mixing condition of a Markov chain such that 
\beq \label{eq:spectral}
\| {\bm W} - {\bf 1}{\bf 1}^\top / N \| \leq 1 - \rho,
\eeq
where $\rho \in (0,1]$ is the spectral gap. Notice that such matrix exists for any connected graph $G$.

With the above mixing matrix, the DGD (or its stochastic variant, DSGD) algorithm \cite{nedic2009distributed,lian2017can} is given by: starting with any $\prm_i^{(0)} \in \RR^d$. $i \in V$, we have
\beq \label{eq:dgd}
\prm_i^{(t+1)} = \sum_{j=1}^N W_{ij} \prm_j^{(t)} - \eta_t {\bm g}_i^{(t)},~\forall~i \in V,~\forall~t \geq 0,
\eeq
where $\eta_t > 0$ is the step size and ${\bm g}_i^{(t)}$ is the local stochastic gradient which is an unbiased estimate of $\grd J_i( \prm_i^{(t)} )$. 
The algorithm \eqref{eq:dgd} operates through a ``consensus-then-optimize'' strategy. At each iteration, the worker $i$ first calculates an average of the local models held by the neighboring workers ($\sum_{j=1}^N W_{ij} \prm_j^{(t)}$), then a stochastic gradient step w.r.t.~the private function ($- \eta_t {\bm g}_i^{(t)}$) is performed. 

We now consider implementing algorithm \eqref{eq:dgd} for decentralized training of an NN model such as \eqref{eq:fi}. As the latter is parameterized by the weights $\prm \in \RR^{mf}$, the workers are required to transmit $mf$ real numbers to the neighbors on the network/graph at each iteration. For $m \gg 1$, this poses a significant challenge due to the limited communication bandwidth. A natural remedy is to compress information before transmitting. 

However, directly compressing the local parameters $\prm_i^{(t)}$ in \eqref{eq:dgd} can result in a non-convergent algorithm since the compression operation can lead to unrecoverable information loss. A better idea is to exploit smoothness and focus on compressing the differences in the iterates. In particular, \cite{koloskova2019decentralized} proposed the CHOCO-SGD algorithm for communication efficient decentalized training, as summarized in Algorithm~\ref{alg:choco}.

\algsetup{indent=.5em}
\begin{algorithm}[h]
	\caption{CHOCO-SGD Algorithm \cite{koloskova2019decentralized}}\label{alg:choco}
	\begin{algorithmic}[1]
		\STATE {\textbf{INPUT}}: initial weights $\{ \prm_i^{(0)} \}_{i=1}^N$, max.~no.~of iterations $T$, consensus parameter $\gamma \in (0,1)$, step sizes $\{ \eta_t \}_{t \geq 0}$. 
		\STATE Set the auxilliary variables $\hat{\prm}_{i,j}^{(0)} = {\bm 0}$, $j \in {\cal N}_i$, $i \in [N]$.
		\STATE Draw the stopping iteration number  ${\sf T} \sim {\cal U} \{0,...,T\}$.
		\FOR{$t=0,1,...,{\sf T}$}
		\FOR[\hfill \texttt{//Local SGD step//}]{$i=1,...,N$}
		\STATE Compute the local SGD: 
		\[ 
			\prm_i^{(t+\frac{1}{2})} = \prm_i^{(t)} - \eta_t {\bm g}_i^{(t)},
		\] 
		where ${\bm g}_i^{(t)}$ is the stochastic estimate of $\grd J_i( \prm_i^{(t)})$.
		\ENDFOR
		\STATE \label{step:comm} For each worker $i=1,...,N$, broadcast the compressed \emph{difference vector} ${\cal Q}( \prm_i^{(t+\frac{1}{2})} - \hat{\prm}_{i,i}^{(t)} )$ to the neighbors, where ${\cal Q}(\cdot)$ is a compression operator satisfying \eqref{eq:compress}.\vspace{.1cm}
		\FOR[\hfill \texttt{//Combination step//}]{$i=1,...,N$}
		\STATE Update the auxiliary variable:
		\[
			\hat{\prm}_{i,j}^{(t+1)} = \hat{\prm}_{i,j}^{(t)} + {\cal Q}( \prm_j^{(t+\frac{1}{2})} - \hat{\prm}_{j,j}^{(t)} ),~\forall~j \in {\cal N}_i. \vspace{-.3cm}
		\]
		\STATE Update the local NN weights:
		\[ \textstyle 
		\prm_i^{(t+1)} = \prm_i^{(t+\frac{1}{2})} + \gamma \sum_{j \in {\cal N}_i} W_{ij} \{ \hat{\prm}_{i,j}^{(t+1)} - \hat{\prm}_{i,i}^{(t+1)} \}.
		\]
		\ENDFOR
		\ENDFOR
		\STATE {\textbf{OUTPUT}}: trained weights $\{ \prm_i^{({\sf T})} \}_{i=1}^N$.
	\end{algorithmic} 
\end{algorithm}

We observe that step~\ref{step:comm} in the algorithm is the only step of CHOCO-SGD involving peer-to-peer communication between the workers. 
To understand the algorithm better, we observe that when the compression operator is an identity operator, i.e.,  ${\cal Q}(\prm) = \prm$ for any $\prm \in \RR^d$, together with the consensus parameter $\gamma = 1$, the CHOCO-SGD algorithm is reduced into:
\beq 
\prm_i^{(t+1)} = \sum_{j=1}^N W_{ij} \{ \prm_j^{(t)} - \eta_t {\bm g}_j^{(t)} \},
\eeq
which resembles the classical DSGD algorithm with a swapped order for the `optimize' and `consensus' steps. 

In general, the communication step depends on a compression operator ${\cal Q}: \RR^d \rightarrow \RR^d$ which reduces the amount of information transmitted. Furthermore, we are compressing the \emph{difference} between the successive iterates with the local stochastic gradient. 
Common communication compression techniques include gradient quantization and gradient sparsification. 
Broadly, these communication compression methods can be classified into \emph{biased} and \emph{unbiased} operators. For the CHOCO-SGD algorithm, it is assumed that the compression operator is a random operator satisfying
\begin{equation}
\label{eq:compress}
	\mathbb{E}_\Omega \left[\| {\cal Q} ( \prm ;\Omega) - \prm \|^{2}\right] \leq (1 - \delta ) \|\prm\|^{2},~~\forall~\prm \in \RR^d,
\end{equation}
where $\Omega$ is the implicit random state of the compression operator, and $\delta \in (0,1]$ is a parameter characterizing the expected error resulted from the compression.
Intuitively, with the condition \eqref{eq:compress}, the CHOCO-SGD algorithm behaves similarly as the DSGD algorithm as only the differences between successive iterates are compressed. 

Examples of \emph{unbiased} ${\cal Q}(\cdot)$ satisfying \eqref{eq:compress} include \cite{alistarh2017qsgd} which \emph{quantizes} a rescaled vector; \cite{sparse1} which \emph{sparsifies} the vectors by setting a random subset of $d-k$ co-ordinates to zero. Concretely, let $\Omega \subseteq \{1,...,d\}$, $|\Omega| = k$ be the selected random subset, we set the operator as ${\sf rand}_k: \RR^d \to \RR^d$ with:
\begin{equation} \label{eq:randk}
	\big[ {\sf rand}_{k}( \prm ; \Omega) \big]_i = \begin{cases}
	[{\prm}]_i & i \in \Omega, \\ 
	0 & i \notin \Omega,
	\end{cases}
\end{equation}
where $[ \prm ]_i$ denotes the $i$th element of the vector $\prm$. 
We observe that ${\sf rand}_{k}( \prm ; \Omega)$ is a vector with at most $k$ non-zero elements.
Notice that for the ${\sf rand}_k$ compression operator, condition \eqref{eq:compress} is satisfied with $\delta = \frac{k}{d}$ \cite{koloskova2019decentralized}.

Meanwhile, \emph{biased} compression operators ${\cal Q}(\cdot)$ are also widely adopted. For example, the sign compression \cite{bernstein2018signsgd} and ${\sf top}_{k}$ sparsification which retains the top-$k$ coordinates in the $d$-dimensional vector with the highest magnitude \cite{sparse1, shi2019understanding, lin2018deep}.
Notably, the ${\sf }_{k}$ sparsification satisfies \eqref{eq:compress} with $\delta = \frac{k}{d}$ since the latter compressor always yield an error lower than that of the ${\sf rand}_k$ sparsifier.

\subsection{Convergence Guarantee and Its Pitfalls} \vspace{.1cm}
We discuss the convergence guarantees of the CHOCO-SGD algorithm and highlights on its pitfalls in decentralized training of overparameterized NNs.

We describe the assumptions used. First, we assume that each private function is smooth, i.e., the gradient map is Lipschitz continuous:
\begin{Assumption} \label{ass:lip}
For any $i \in [N]$, there exists $L \geq 0$ such that 
\beq
\| \grd J_i( \prm ) - \grd J_i( \prm' ) \| \leq L \| \prm - \prm' \|,~\forall~\prm, \prm' \in \RR^d. 
\eeq
\end{Assumption}
Next, we specify conditions on the stochastic gradient estimates. Denote ${\cal F}_t$ as the filtration generated by the random variables $\{ \prm_i^{(\tau)} : i \in V, 0 \leq \tau \leq t \}$. We assume that:
\begin{Assumption} \label{ass:stoc}
There exists $\sigma, G \geq 0$ such that for any $i \in [N]$, $t \geq 0$, the stochastic gradient ${\bm g}_i^{(t)}$ satisfies
\beq 
\begin{split}
& \EE[ {\bm g}_i^{(t)} | {\cal F}_t ] = \grd J_i( \prm_i^{(t)} ),~~\EE[ \| {\bm g}_i^{(t)} \|^2 | {\cal F}_t ] \leq G^2, \\
& \EE[ \| {\bm g}_i^{(t)} - \grd J_i( \prm_i^{(t)} ) \|^2 | {\cal F}_t  ] \leq \sigma^2.
\end{split}
\eeq
\end{Assumption}
Notice that the first condition is satisfied by a uniform sampler for the (mini-batch) stochastic gradient, e.g., when ${\bm g}_i^{(t)} = \grd_{\prm} {\sf loss}( f( {\bm x}_{j_t}; \prm_i^{(t)}) ; y_{j_t} )$ such that $j_t$ is selected uniformly at random from $\{ 1, ..., |{\cal D}_i| \}$. We have also assumed that the stochastic gradients have bounded second order moments.

\begin{figure*}[th]
	\vspace{-0.5em}
	\centering
	\subfigure[
		${\sf top}_k$ sparsification.
	]{
		\includegraphics[width=.45\textwidth,]{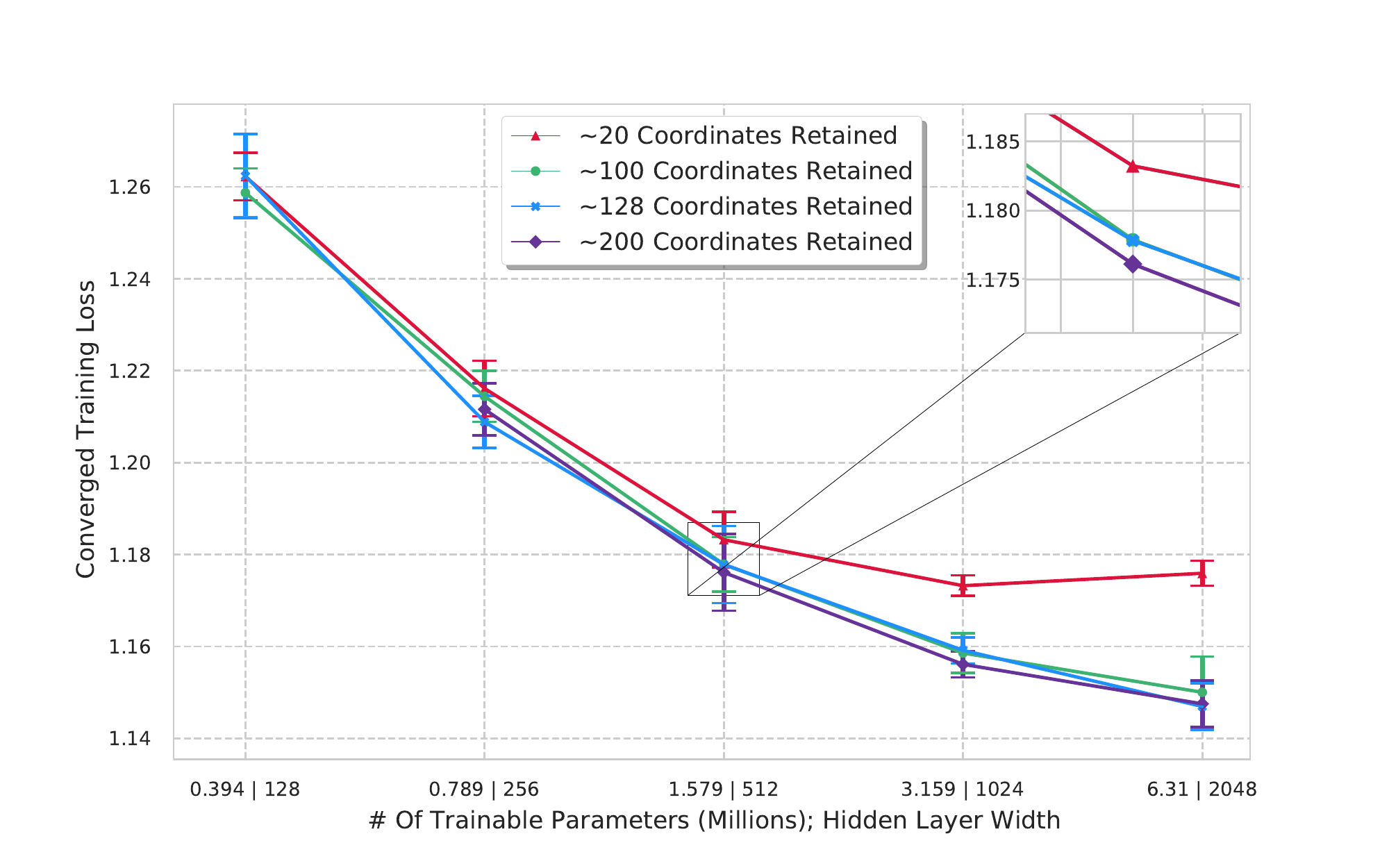}
		\label{fig_topk}
	}
	\hfill
	\subfigure[
		${\sf rand}_k$ sparsification
	]{
		\includegraphics[width=.45\textwidth,]{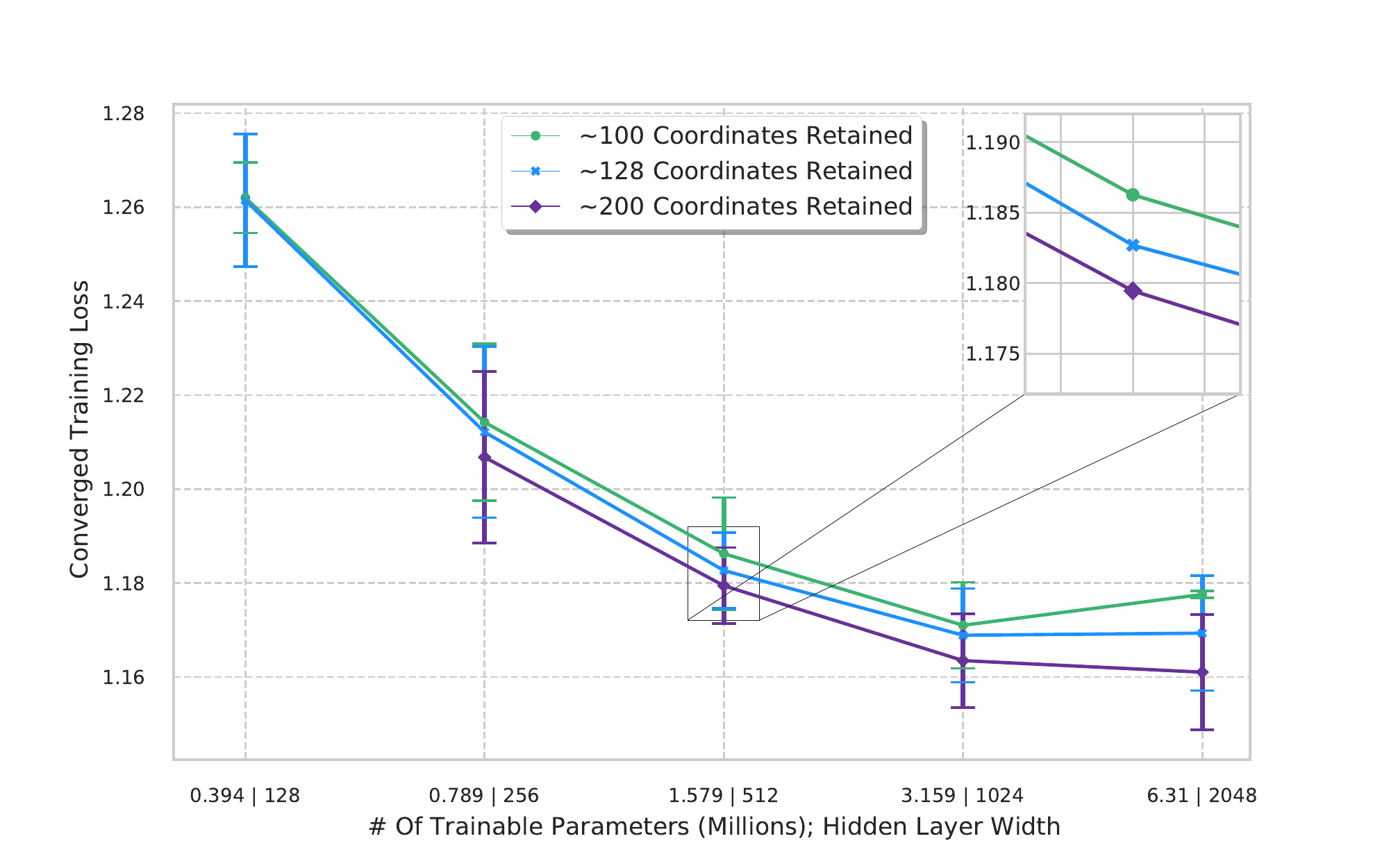}
		\label{fig_randk}
	}
	\hfill
	\vspace{-0.4em}
	\caption{
		Converged training losses for large-width NNs compressed with biased ${\sf top}_k$ \ref{fig_topk} and unbiased ${\sf rand}_k$ \ref{fig_randk} sparsification on CIFAR-10. Reported accuracies averaged over $5$ independent trials with error bars indicating $99\%$ confidence interval over the $5$ trials. For $k =20$ with ${\sf rand}_k$ sparsification, divergence is observed in the case of 2048-layer and 1024-layer NNs on all $5$ trials.
	} \label{fig:cifar}
\end{figure*}

We observe the following result that is borrowed from \cite[Theorem 4.1]{Koloskova*2020Decentralized} on the CHOCO-SGD algorithm:
\begin{Theorem} \label{th:main}
Under Assumptions~\ref{ass:lip}, \ref{ass:stoc} and suppose that the compressor satisfies \eqref{eq:compress}. There exists $\eta, \gamma > 0$ such that if we consider a constant step size with $\eta_t \equiv \eta$, then for any $T \geq 1$, the output generated by Algorithm~\ref{alg:choco} satisfy:
\beq \notag
\EE[ \| \grd J( \overline{\prm}^{({\sf T})} ) \|^2 ] = {\cal O} \left( \sqrt{\frac{L \sigma^2 J_0}{N T} } + \left( \frac{L G J_0}{\rho^2 \delta T } \right)^{\frac{2}{3}} \right),
\eeq
where the expectation is taken over ${\sf T}$ and the stochastic quantities in the algorithm, $\delta$ was defined in \eqref{eq:compress}, $\rho \in (0,1]$ is the spectral gap of ${\bm W}$ defined in \eqref{eq:spectral}, $\overline{\prm}^{({\sf t})} = \sum_{i=1}^N \prm_i^{({\sf t})} / N$ is the network average iterate, and $J_0 = J( \overline{\prm}^{(0)} ) - \min_{\prm} J(\prm)$.
\end{Theorem}
The theorem suggests that the CHOCO-SGD algorithm finds an ${\cal O}(1/\sqrt{T})$-stationary solution to \eqref{eq:opt_int} in at most $T$ iterations as we note that ${\sf T} \leq T$.\vspace{.2cm}

\noindent \textbf{Convergence for Overparameterized Models.}~~
We concentrate on the performance of CHOCO-SGD when $d \gg 1$, for example, when training an overparameterized NN such as \eqref{eq:fi} with $m \gg 1$ neurons. Furthermore, to control the bandwidth usage, we choose the ${\sf rand}_k$ sparsifier or ${\sf top}_k$ sparsifier as the compressor. Now, we fix the number of iterations as $T$, and the number of coordinates sent per iteration at $k$, i.e., we fix the amount of data transmitted in the CHOCO-SGD algorithm. In this setting, we have 
\beq \label{eq:thm1d}
\EE[ \| \grd J( \overline{\prm}^{({\sf T})} ) \|^2 ] = 
{\cal O} \Big( \sqrt{\frac{L \sigma^2 J_0}{N T} } + d^{\frac{2}{3}} \Big( \frac{LG J_0}{\rho^2 k T } \Big)^{\frac{2}{3}} \Big) .
\eeq
Furthermore, applying Theorem~\ref{th:main} shows that to reach an $\epsilon$-stationary solution (i.e., $\EE[ \| \grd J(\overline{\prm}^{({\sf T})} ) \|^2 ] \leq \epsilon$), the number of CHOCO-SGD iterations required grows in the order: 
\beq \label{eq:pitfall}
T = \Omega \left( L J_0 \cdot \max\left\{ \frac{ \sigma^2 }{ N \epsilon^2 } \, , \,  \frac{d}{k} \frac{ G }{ \rho^2 \epsilon^{1.5} } \right\} \right)  .
\eeq
As the amount of data transmitted per iteration is constant (i.e., $k$ real numbers), the above calculation indicates that \emph{the CHOCO-SGD algorithm may require a higher communication complexity as the NN model becomes inncreasingly overparameterized (i.e., when $d \gg 1$)}, in order to maintain the same performance level, despite the compression being applied at each iteration. In fact, for any $1 \leq k \leq n$, it is predicted from \eqref{eq:pitfall} that the number of real numbers transmitted is $\Omega( \max\{ k \sigma^2 / (N \epsilon^2), d G / (\rho^2 \epsilon^{3/2}) \} )$ when the ${\sf top}_k$ or ${\sf rand}_k$ sparsifier is used as the compressor in CHOCO-SGD. We notice that similar dependence on the problem dimension $d$ is also observed in other compressed DSGD methods, e.g., \cite{kovalev2021linearly,tang2018communication,alistarh2017qsgd}. 

The observations in \eqref{eq:thm1d}, \eqref{eq:pitfall} indicate a pitfall in the existing theories when applying compressed DSGD methods such as CHOCO-SGD to overparameterized NNs. Particularly, it suggests that although the compression scheme can reduce the communication cost \emph{per iteration}, the number of iterations required would be increased. 
In the next section, we conduct an extensive set of numerical experiments to test the above observation. Interestingly, we show that in most cases, \emph{increasing the degree of overparameterization can lead to a better performance for the trained NN without increasing the communication cost during training}, which is in contrast to the said observation.
We conjecture that such phenomena is a consequence of an artifact in the existing analysis and discuss the possible fix for the observation.

\section{Empirical Studies}\label{sec:num}
In this section, we perform numerical experiments to examine the performance of compressed DSGD method when applied to training overparameterized NNs. We concentrate on the effects of the number of parameters on the communication complexity using the CHOCO-SGD method \cite{Koloskova*2020Decentralized} [cf.~Algorithm~\ref{alg:choco}]. 
For simplicity, we consider a two-layer NN with ReLU activation described in \eqref{eq:fi} and adjust the width, $m$, of the NN. 

We consider the task of training a classifier with the CIFAR-10 dataset \cite{krizhevsky2009learning} that contains 50K (resp.~10K) training (resp.~test) samples. Each sample consists of a $32 \times 32$ RGB image which can be represented as a $3072$-dimensional vector, and is associated with a label selected from 10 image classes. To simulate the decentralized training environment, samples from the 10 image classes are uniformly split among $N$ workers and shuffled at every epoch – as in \cite{Koloskova*2020Decentralized, goyal2017accurate}. 
To establish a challenging generalization task, we test the trained models on CIFAR-10.1 \cite{recht2018cifar}. 

For the CHOCO-SGD method, we use a minibatch size of $\xi = 128$ for every iteration. 
Our chosen mode of communication compression is ${\sf top}_k$ and ${\sf random}_k$ with a fixed number of co-ordinates $k$ allowed to be communicated between workers. We choose a constant consensus parameter $\gamma = 0.0375$ in Algorithm~\ref{alg:choco} and an SGD stepsize $\eta = 0.1$ which is decreased by a factor of $10$ on epochs $100,150, 200$. Our ${\sf top}_k$ and ${\sf random}_k$ simulations are run on $N=8$ nodes of a ring topology. The decentralized training environment is simulated on an MPI-based \cite{graham2005open} network where we assign an independent CPU process to each worker. 

\begin{figure}[t]
	\vspace{-0.1em}
	\centering
		\includegraphics[width=.45\textwidth,]{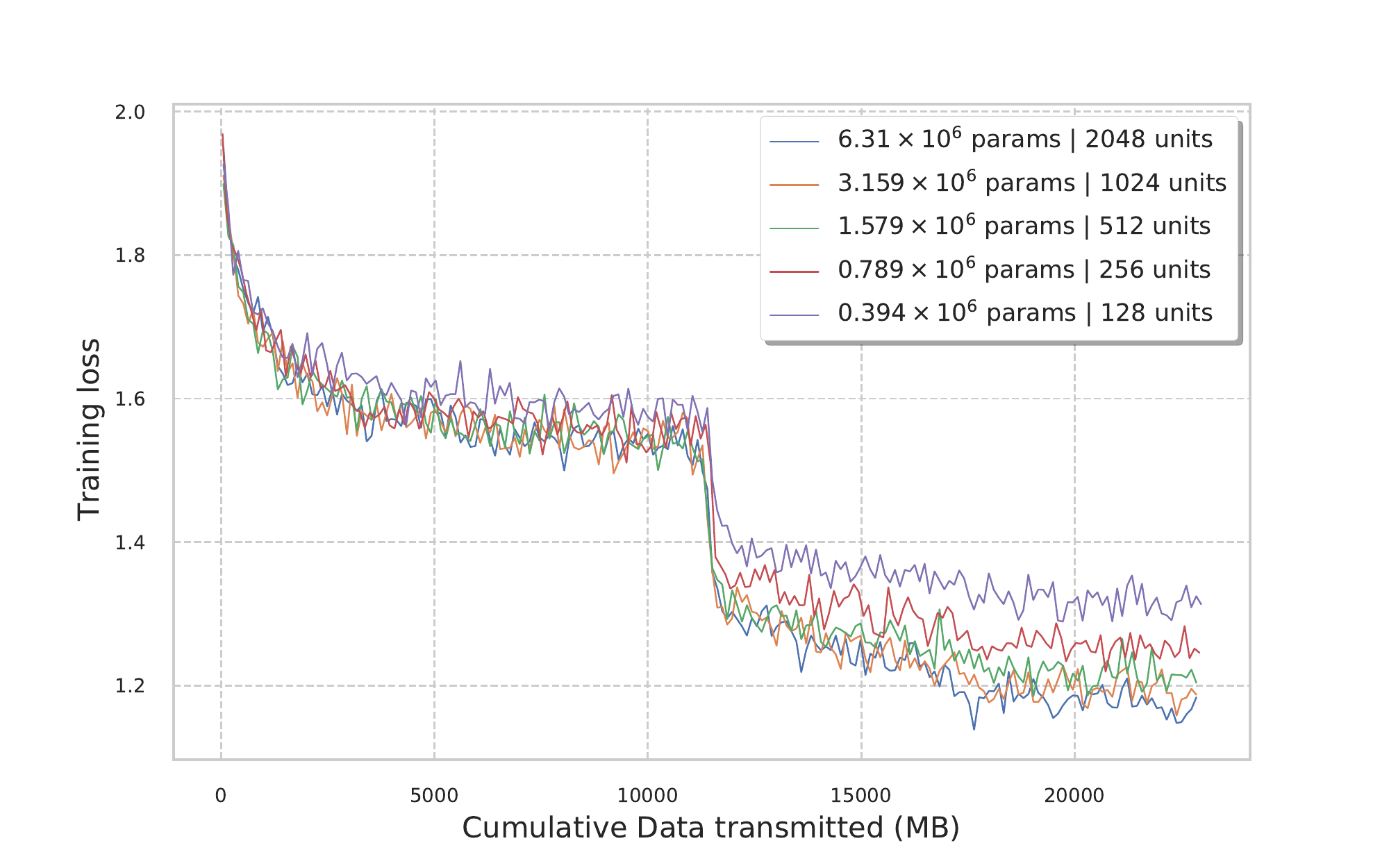}
	\vspace{-0.4em}
	\caption{
		Training loss with cumulative communication cost in (MB) for large-width NNs of varying widths constrained with a constant sparsification coordinate $k = 100$. While \emph{convergence rate} of all models is invariant to layer width, overparameterized models ($d \geq 2 \times 10^6$ parameters) converge to a solution of lower training loss with identical cumulative data usage.
	} \label{fig:data}
\end{figure}

Our first empirical example aims at comparing the quality of solution found by CHOCO-SGD after performing {$300$} epochs of iterations against the number of parameters of the trained NN ($d$). Notice that as $k$ is fixed in the compressor design, the communication complexity (i.e., number of bits transmitted over the network) is fixed. 
Our results are presented in Fig.~\ref{fig_topk} and \ref{fig_randk}. From the figures, we observe that the training loss \emph{decreases} with $d$, especially when the ${\sf top}_k$ sparsifier is used. 
Notice that the above is actually \emph{in contrast} to the observations made in \eqref{eq:thm1d} based on Theorem~\ref{th:main}, which predicts that the solution quality decreases as the problem dimension increases. 
In Table~\ref{tab:cifar} we compare the \emph{converged} test accuracy evaluated on CIFAR10.1 with different widths for the NN, where similar observations are obtained as in Fig.~\ref{fig:cifar}. 
Lastly, we plot the trajectory of the training loss against the cumulative communication cost (in MB) in Fig.~\ref{fig:data}. Again, we observe that increasing model dimension (overparameterization) leads to faster convergence with respect to the communication cost. 

\input{figures/table_sparse.tex}



\begin{figure}[t]
	\centering
		\includegraphics[width=.45\textwidth,]{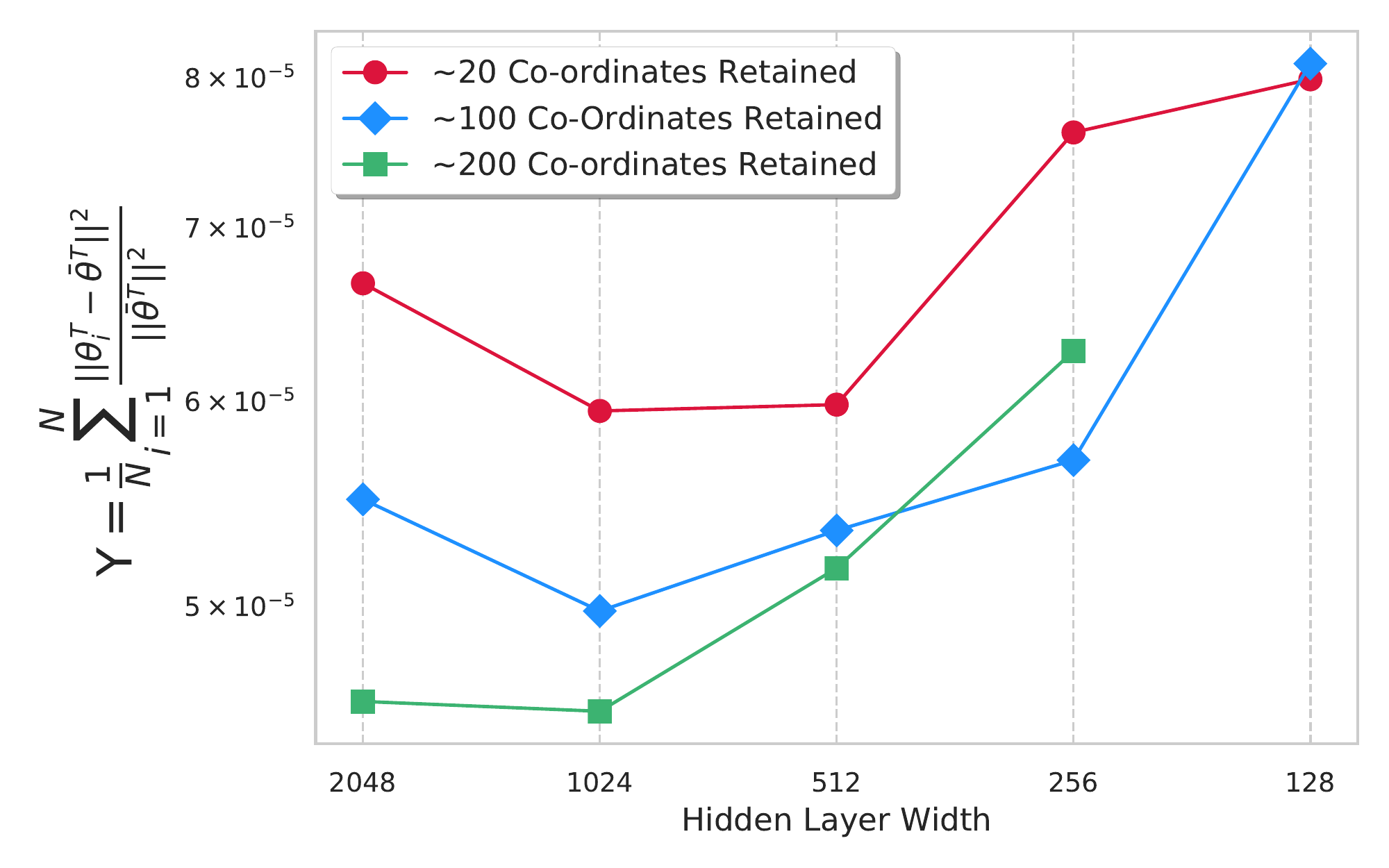}
	\vspace{-0.4em}
	\caption{
		Converged, normalized consensus distance between $N=8$ workers for different NN layer widths. Overparameterized models enjoy significantly greater consensus among workers with only marginal dependence on sparsification co-ordinate bandwidth for larger models. 
	} \label{fig:consensus_normalized}
\vspace{-0.5em}
\input{figures/consensus_rate}
\vspace{-2em}
\end{figure}

Our second empirical example examines the consensus error in the converged solution after {$300$} epochs of CHOCO-SGD with $N=8$ workers. We aim at studying the consensus error of the converged solutions and the effects of problem dimension. Here, the normalized consensus error is defined as:
\beq \label{eq:cons_dist}
\Upsilon = \frac{1}{N}\sum_{i=1}^{N} \frac{ \|\prm_{i}^T - \overline{\prm}^T \|^{2} } { \| \overline{\prm}^T \|^2 } ,
\eeq
where $\prm_{i}^T$ denotes the CHOCO-SGD solution after 300 epochs.  
The above metric is compared against the width of the hidden layer, $m$, in the NN in Fig.~\ref{fig:consensus_normalized}.  
Naturally, we observe that the consensus error increases when the number of coordinates kept in the sparsifier $k$ decreases. 
Moreover, an intriguing observation is that the consensus error \emph{increases} when the width of the NN \emph{increases} from $m=1024$ to $m=2048$. We recall from Fig.~\ref{fig:cifar}, Table~\ref{tab:cifar} that the training loss/testing accuracy actually \emph{decreases} with this change in the number of parameters. This suggests that a consensual solution may not be necessary for achieving a lower training loss in the overparameterized setting. 



\begin{figure}[t]
	\vspace{-0.4em}
	\centering
		\includegraphics[width=.475\textwidth]{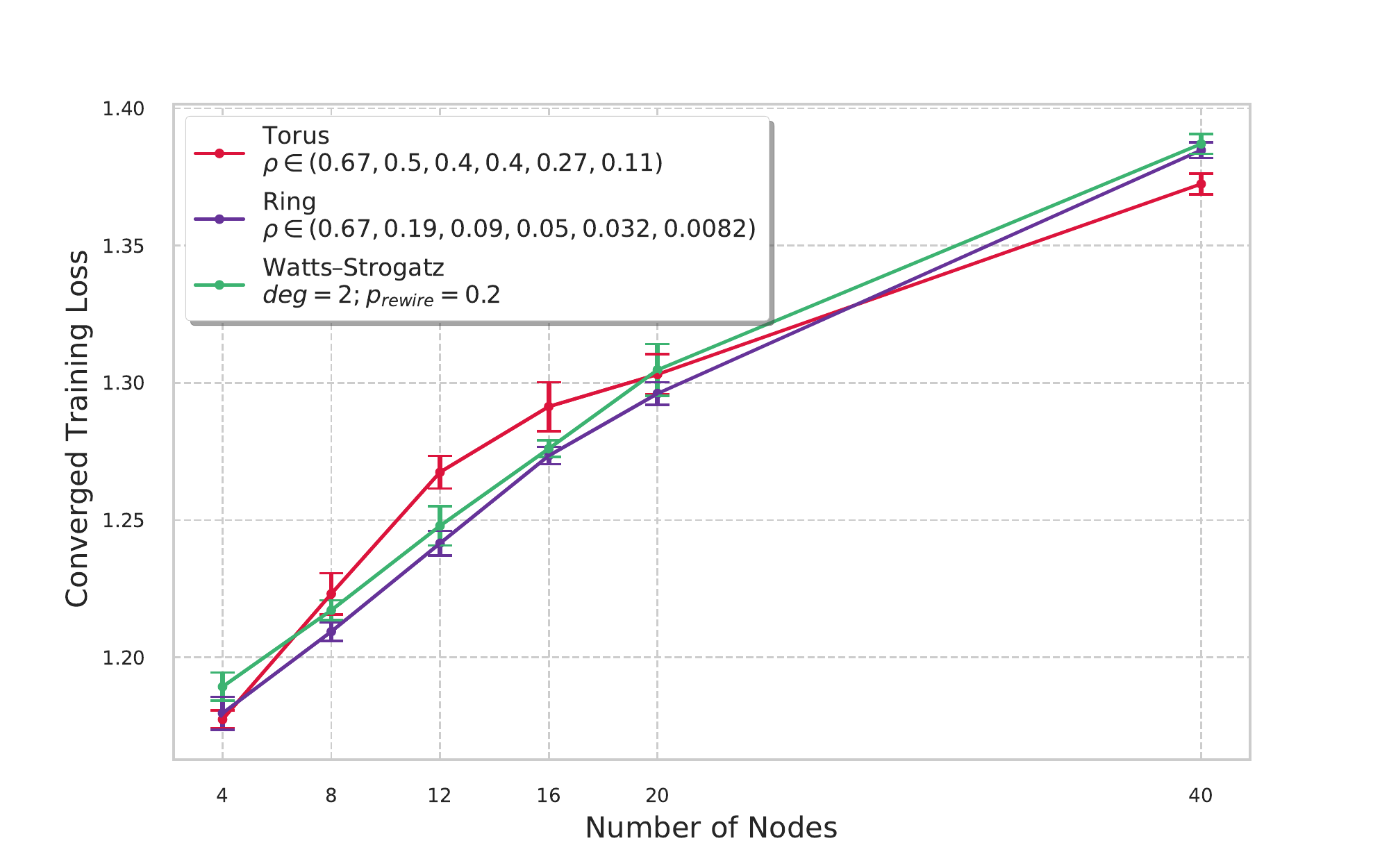}
	\vspace{-0.6em}
	\caption{
		Converged training loss versus number of workers with different graph topologies. 
		Error bars indicate 95\% confidence interval over $5$ independent trials. 
	} \label{fig:topology}
\end{figure}

Our third empirical example examines the effects of graph topology on the converged training loss after {300} epochs of CHOCO-SGD iterations. In this experiment, we fix the number of hidden neurons at {$m= 256$ and use a ${\sf top}_k$ sparsifier with $k = 51$}. We experiment with common graph topologies such as ring $(deg_{ring} = 2)$, torus $(deg_{torus} = 4)$, and Watts–Strogatz (small-world) graphs \cite{watts1998collective}, and for different numbers of connected neighbors. Notice that the spectral gap parameter $\rho$ decreases with the number of workers $N$, especially for the ring topology. The result is presented against the number of workers in Fig.~\ref{fig:topology}. As observed in the figure, for all topology settings, the training loss increases with the number of workers $N$. The deterioration in performance is the most severe with the ring topology. 
Unlike the previous example, we notice that the observed behavior is consistent with the theory in Theorem~\ref{th:main}, which predicts a slower convergence when the graph topology has a smaller $\rho$.

\begin{figure}[t]
	\vspace{-0.4em}
	\centering
		\includegraphics[width=.475\textwidth]{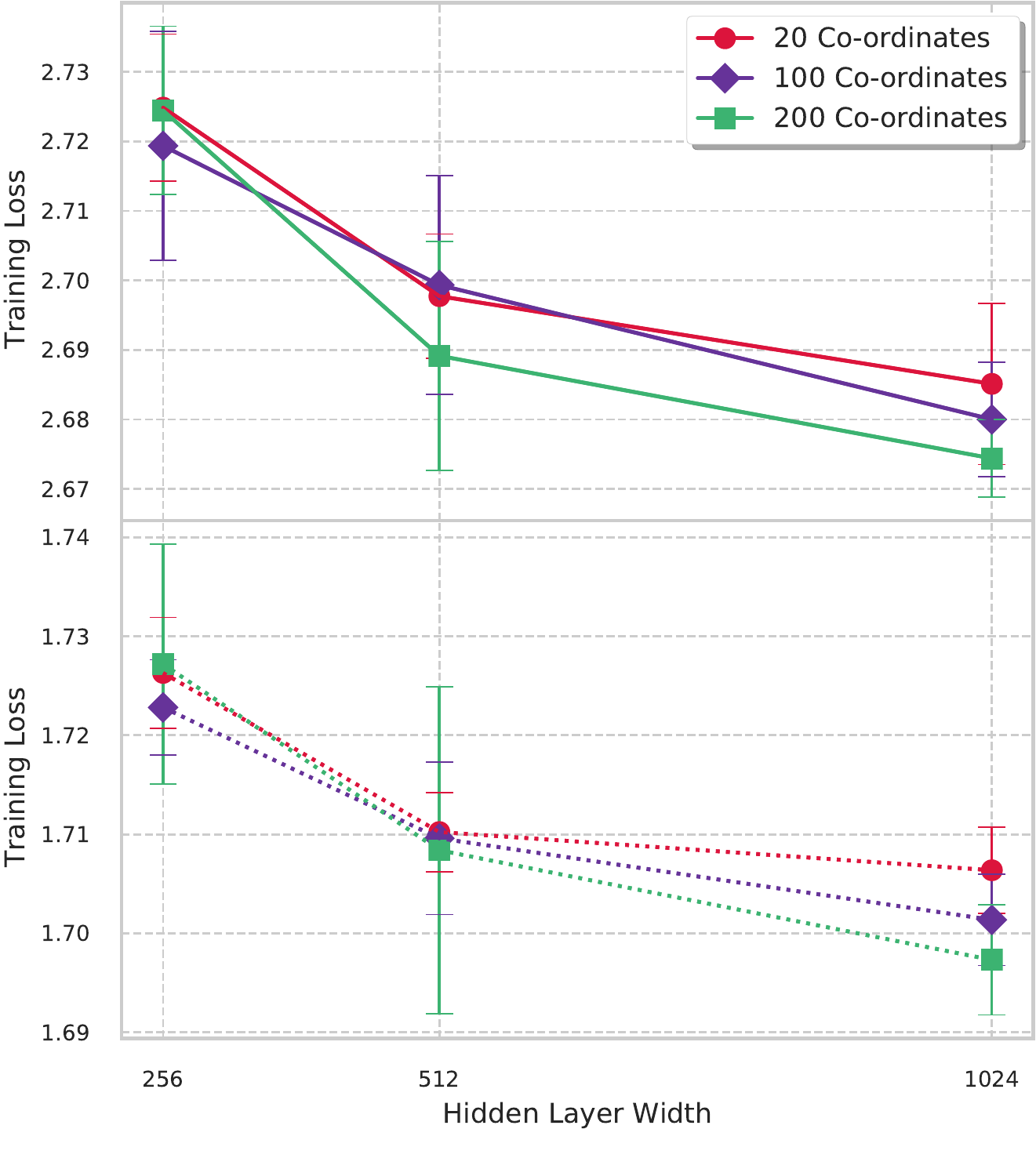}
	\vspace{-0.6em}
	\caption{
		Converged training loss versus model dimensionality for \emph{heterogeneously} distributed data among $N=8$ workers on a ring topology where each worker's private distribution contains labels from a single class. Top row (solid line) indicates average performance of the \emph{local model} on the complete training dataset. Bottom row (dotted line) indicates performance of \emph{averaged model} ($\bar{\theta}^{*}$) on the complete training dataset. Error bars indicate $95\%$ confidence interval over 5 independent trials. 
	} \label{fig:iid}
\end{figure}

Lastly, we study a case of decentralized training with heterogeneous data. Particularly, we allocate samples from a single class to each of the $N=8$ workers. We evaluate the training loss with \emph{global data} using the local model $\prm_i^T$ obtained by CHOCO-SGD with ${\sf top}_k$ sparsifier after 300 epochs; compared to the averaged model $(1/N)\sum_{i=1}^N \prm_i^T$. The simulation result is shown in Fig.~\ref{fig:iid}. We observe that a similar trend as Fig.~\ref{fig:cifar} holds as the training loss decreases with the hidden layer width $m$ given that the algorithms are run with a similar communication budget. Furthermore, we note that the training loss is higher with the local model.
\vspace{.2cm}



\noindent \textbf{Discussions.}~~ The above numerical examples demonstrate a consistent discrepancy between existing theories on compressed DSGD methods and their practical performances when training overparameterized NN models. In particular, we show that with the same communication cost allowed, the performance of the CHOCO-SGD trained NN improves with the number of neurons $m$ employed in the NN, contrary to \eqref{eq:thm1d}.

The careful readers may notice that the poor dependence on dimension $d$ in the convergence bounds \eqref{eq:thm1d} can be seen as a direct consequence of the worst-case analysis of the compression error in \eqref{eq:compress}. This is because one is forced to take $\delta = d / k$ to account for the condition that holds for all ${\bm x} \in \RR^d$ in case of the ${\sf top}_k$/${\sf rand}_k$ sparsifier. Furthermore, we remark that the existing theories are developed for general distributed optimization problems, assuming only high level properties such as smoothness of the objective function. 

On the other hand, recent works on overparameterized models \cite{jacot2018neural, arora2019fine, bakshi2019learning} have shown that the training of such models with SGD shall be compared to the learning of a kernel model in the RKHS. Among others, a key innovation therein is that the optimality gap is measured in terms of the distance in the (infinite-dimensional) function space, where the rate of convergence is shown to be independent of $m$ even as $m \rightarrow \infty$. We believe that extending these results to a decentralized setting such as the CHOCO-SGD method can help break the curse of dimensionality in the analysis. 
We remark that an interesting observation was made in \cite{shevchenko2020landscape} on the stability of overparameterized NN models with dropout, an operation that is akin to applying the sparsifier. 

\section{Conclusions}
This paper provides the first empirical study on the performance of compressed DSGD methods on \emph{overparameterized NN models}. On the positive side, our result shows that utilizing overparameterized NN models in a decentralized learning is both practical and beneficial, contradicting existing theories that suggest otherwise. Furthermore, we identify a gap in the existing theories that attempt to analyze compressed DSGD methods. We believe that our article will serve as a motivating study to develop improved communication efficient algorithms for decentralized training of NNs.  

\section*{Acknowledgements}
This work is supported by CUHK Direct Grant \#4055113.
\bibliographystyle{IEEEtran}
\bibliography{refs}
\end{document}

%% file: figures/table_sparse.tex
\begin{table*}[!t]
	\centering
	\small
	
	\resizebox{1\textwidth}{!}{%
		\huge
		\begin{tabular}{ccccc|ccc}
			\toprule
			\multirow{2}{*}{\diagbox{\# Neurons}{\# Coordinates $k$}}  & \multicolumn{4}{c|}{${\sf top}_k$ sparsifier}            & \multicolumn{3}{c}{${\sf rand}_k$ sparsifier}  \\  \cmidrule(lr){2-8}
			       & 20 & 100 & 128 & 200   & 100 & 128  & 200     \\ \hline
			$m=2048$ &      46.481 $\pm$ 0.767 & 
							46.837 $\pm$ 0.354 &  
							46.663 $\pm$ 1.544 &  
							46.711 $\pm$ 0.755 & 
							
							46.919 $\pm$ 1.401 & 
							43.288 $\pm$ 0.769 & 
							46.365 $\pm$ 1.1 \\

			$m=1024$ & 		46.776 $\pm$ 0.423 & %
							46.184 $\pm$ 0.937 & %
							46.688 $\pm$ 1.05 & %
							47.156 $\pm$ 0.999 & 
							
							48.245 $\pm$ 0.855& 
							47.135 $\pm$ 1.103 & 
							47.309 $\pm$ 0.984 \\
			
			$m=512$ & 		45.651 $\pm$ 2.219 & 
							46.266 $\pm$ 0.902 & 
							45.787 $\pm$ 1.111 & 
							46.292 $\pm$ 0.584 & 
							
						    46.497 $\pm$ 0.653& 
							46.47 $\pm$ 0.401 & 
							45.505 $\pm$ 0.81 \\
									
			$n=256$ & 		45.286 $\pm$ 0.413 & 
							45.148 $\pm$ 1.103 & 
							45.051 $\pm$ 0.802 & 
							45.123 $\pm$ 0.855 & 
							
							45.707 $\pm$ 0.751 & 
							44.928 $\pm$ 0.846&
							44.733 $\pm$ 0.616  \\

			$m=128$ & 		44.59 $\pm$ 0.984 & 
							43.967 $\pm$ 0.851 & 
							42.778 $\pm$ 0.45 & %
							N/A &
							
							43.288 $\pm$ 1.435 & 
							43.568 $\pm$ 1.282 & 
							N/A \\
			\hline 
			comm.~cost~per~iteration~(MB)& 0.482 & 2.410 & 3.012 & 4.820 & 2.410 & 3.012 & 4.820 \\
			\bottomrule \\
		\end{tabular}%
		\vspace{-2em}
	}
	\caption{\textbf{Generalization performance of CHOCO-SGD to CIFAR-10.1 with overparameterized models}. Contrary to convergence performance, large-width NNs trained with ${\sf top}_k$ and ${\sf rand}_k$ sparsification exhibit marginal dependance to model dimension at test-time (average accuracy disparity of $2.39\%$ between $2048$ and $128$ unit models). Test results averaged over $5$ independent trials with altered seeds on a ring-topology communication graph with $8$ nodes. CHOCO-SGD diverges when parameters are compressed with ${\sf rand}_k$ sparsification of $k = 20$ co-ordinates.} \label{tab:cifar}

\end{table*}

%% file: figures/consensus_rate.tex
\begin{table}[H]
\centering 
\vspace{1.2em}
\resizebox{\columnwidth}{!}{%
 \begin{tabular}{ccccc}
        \toprule
    \multirow{2}{*}[-3pt]{Layer Width}  & \multicolumn{3}{c}{Normalized Consensus Distance [cf.~\eqref{eq:cons_dist}]} \\ 
        \cmidrule{2-5} 
      & Epoch $=200$ & Epoch $=100$ & Epoch $=50$ & \\ 
         \midrule
    2048  & $5.499 \times 10^{-5}$  & $9.8206 \times 10^{-3}$ & $1.3977 \times 10^{-2}$ \\ 

    1024  & $4.980 \times 10^{-5}$ & $1.0346 \times 10^{-2}$ & $1.5307 \times 10^{-2}$ \\ 

    512  & $5.349 \times 10^{-5}$ & $1.0026 \times 10^{-3}$ & $1.3478 \times 10^{-2}$ \\ 
    
    256  & $5.694 \times 10^{-5}$ & $8.7639 \times 10^{-3}$ & $1.2423 \times 10^{-2}$ \\
    
    128  & $8.098 \times 10^{-5}$ & $7.3181 \times 10^{-3}$ & $9.2698 \times 10^{-3}$ \\
        \bottomrule
    \end{tabular} 
    }
\caption{\footnotesize \textbf{Converged consensus distances at intermediate training epoch numbers.} Overparameterized models converge to better-consensus solutions at a slower rate compared to low-width NNs}
\end{table}